\pdfoutput=1
\documentclass{article}
\usepackage[margin=25mm]{geometry}
\usepackage{amsmath}
\usepackage{amsfonts}
\usepackage{amssymb}
\usepackage{graphicx}

\usepackage[utf8]{inputenc} 
\usepackage[T2A]{fontenc}
\usepackage[english]{babel}
\usepackage[unicode, pdftex]{hyperref}

\usepackage{xcolor}
\usepackage[all]{xy}

\newtheorem{theorem}{Theorem}

\title{Planar trees as complete topological invariants of Morse flows with a sink on the 2-sphere }
\author{Oleksandr Pryshliak}

\date{\today}

\begin{document}
\maketitle

\begin{abstract}
To investigate the topological structure of Morse flows with a sink on the 2-sphere we use the planar tree as complete topological invariant of the flow. 
We give a list of all planar tree with at least 7 edges. We use a list of rooted planar trees.
\end{abstract}

\subsection*{Introduction}
\subsection{Relevance of the topic}

The study of the topological structure of embedded trees on a plane is a topical research topic for several reasons:

Applications in various scientific fields: Topological tree embeddings on a plane have a wide range of applications in various scientific fields, such as computer graphics, computer visualization, bioinformatics, social networks, transportation networks, and many others. The study of topological structures that arise when nesting trees on a plane can help to understand the interactions between objects in different contexts and use this knowledge in various scientific and practical applications.

Challenges and problems: Embedded trees on a plane is a non-trivial problem, as it requires taking into account a number of restrictions, such as restrictions on the intersection of edges, restrictions on the distance between nodes, and others. Understanding the topological structure of nesting trees on a plane can help solve complex problems and challenges related to the development of efficient nesting algorithms, optimizing interactions between objects, and improving data visualization.

Development of scientific methods: The study of the topological structure of nesting trees on a plane contributes to the development of scientific methods and approaches, such as graph theory, computer geometry, optimization methods, and others.

\subsection{Literature review}

In the \cite{fra90} definition
basic concepts and properties of topological nesting of trees on a plane, such as the concepts of "adjacency", "levels", "visibility" and "external face". The authors developed the concept of "canonical embeddings" and established a connection between them and the topological properties of trees.

In \cite{che91}, the problem of placing the vertices of trees on a plane in such a way that the edges are straight segments was investigated. The concept of "straight-line visible edges" was introduced and it was shown that not every tree has direct attachments on the plane.

In \cite{mil91}, the authors studied the structural properties of nesting trees on a plane, in particular, for the first time defined the concept of "planar separators" (planar separators) and showed how to use them to construct parallel nestings of trees on a plane.

  In \cite{buc02}, the authors proposed a method for constructing "beautiful" nestings of trees on a plane.

Book\cite{bat99}
  provides a general overview of algorithms for graph visualization, including techniques for placing tree vertices on a plane.

\cite{tam13} can be considered a reference publication that contains articles from leading researchers in the field of graph drawing and visualization. In it, you can find sections devoted to the topological properties of nesting trees on the plane and the corresponding algorithms.

In the collection of papers \cite{ari17} presented at the Graph Drawing 2016 conference, the authors presented the latest research in the field of graph drawing, including work on the topological structure of nesting trees on a plane.

Work \cite{pac01}
is devoted to theory and algorithms for working with planar graphs, including placing tree vertices on the plane.

In \cite{kob13}
the author offers various methods and models for placing graph vertices on a plane.

The application of nested trees to study is especially interesting
topological properties of functions and dynamic systems -- such properties that do not depend on the choice of coordinate systems, measurement methods, etc. For typical systems, they do not change with small changes in the parameters defining this system.

Topological classification of flows on closed surfaces was obtained in  \cite{bilun2023gradient, Kybalko2018, Oshemkov1998, Peixoto1973, prishlyak1997graphs, prishlyak2020three, akchurin2022three, prishlyak2022topological, prishlyak2017morse, kkp2013, prishlyak2021flows, prishlyak2020topology, prishlyak2019optimal, prishlyak2022Boy},
  and on surfaces with a boundary y
\cite{bilun2023discrete, bilun2023typical, loseva2016topology, loseva2022topological, prishlyak2017morse, prishlyak2022topological, prishlyak2003sum, prishlyak2003topological, prishlyak1997graphs, prishlyak2019optimal, stas2023structures}.
Complete topological invariants of Morse-Smale flows on 3-manifolds was constructed in \cite{prish1998vek, prish2001top, Prishlyak2002beh2, prishlyak2002ms, prishlyak2007complete, hatamian2020heegaard, bilun2022morse, bilun2022visualization}.

Morse streams are gradient streams for Morse functions. If we fix the value of the function at special points, then the structure of the flow determines the structure of the function \cite{lychak2009morse, Smale1961}. Therefore, Morse--Smale flows classification is related to the classification of the Morse functions.

Topological invariants of functions on oriented surfaces were constructed in \cite{Kronrod1950} and \cite{Reeb1946} and in \cite{lychak2009morse} for unoriented surfaces, and in \cite{Bolsinov2004, hladysh2017topology, hladysh2019simple, prishlyak2012topological} for surfaces with a boundary, in \cite{prishlyak2002morse} for non-compact surfaces.

Topological invariants of smooth functions were also studied in \cite{bilun2023morseRP2, bilun2023morse, hladysh2019simple, hladysh2017topology, prishlyak2002morse, prishlyak2000conjugacy, prishlyak2007classification, lychak2009morse, prishlyak2002ms, prish2015top, prish1998sopr, bilun2002closed, Sharko1993}, for manifolds with limit in \cite{hladysh2016functions, hladysh2019simple, hladysh2020deformations}, and on 3- and 4-dimensional manifolds in \cite{prishlyak1999equivalence, prishlyak2001conjugacy}.

For a first introduction to the topological theory of functions and dynamical systems, we recommend  \cite{prishlyak2012topological, prish2002theory, prish2004difgeom, prish2006osnovy, prish2015top}.

\subsection{The purpose of the paper}

The main goal of this work is to describe all possible non-isomorphic planar trees with no more than 8 vertices and use it for topological classifications of Morse flows with one sink on $S^2$.

\subsection{Research methods}

To solve the main problem, we will use known algorithms for the classification of planar-rooted trees. We will call the center of the tree the vertex that has the greatest distance to the vertices of valence 1. The center of the tree divides it into root trees. Our research method consists of gluing the roots of root trees.

\section{Planar-rooted trees}

Plane-rooted trees are a particular type of tree in which the vertices are located on a plane in such a way that the edges do not intersect. In these trees, a special vertex is defined, known as the root of the tree, from which all the edges emanate. Each vertex can have several child vertices that are located lower on the plane than their parent vertices.

\begin{figure}[ht!]
\center{\includegraphics[width=0.2\linewidth]{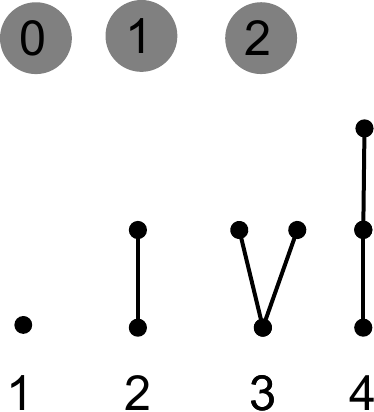}
}
\caption{rooted trees with 0, 1, 2  and 3 edges
}
\label{t0-2}
\end{figure}

Planar-rooted trees are important data structures used in many fields such as computer science, graphic design, data visualization, data processing algorithms, and many others. They are one of the basic data structures for representing a hierarchical organization of data, where one element (the root) has a special status relative to other elements.

\subsection{Rooted trees with 0, 1 and 2 edges
}
We will use the list of root trees given in the work \cite{stas2023structures}. In fig. \ref{t0-2} lists all possible rooted trees with no more than 2 edges (1+1+2 trees).

\subsection{Rooted trees with 3 edges}

\begin{figure}[ht!]
\center{\includegraphics[width=0.27\linewidth]{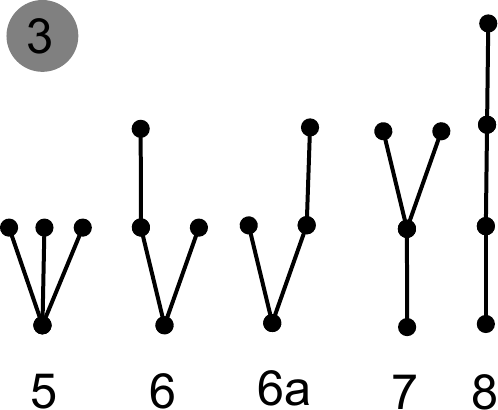}
}
\caption{rooted trees with 3 edges
}
\label{t3}
\end{figure}
In fig. \ref{t3} all possible 5 rooted trees with 32 edges are given. Graphs 6 and 6a are isomorphic on the plane, but not isomorphic on the half-plane (the root is located on the boundary of the half-plane).

\subsection{Rooted trees with 4 edges}
\begin{figure}[ht!]
\center{\includegraphics[width=0.25\linewidth]{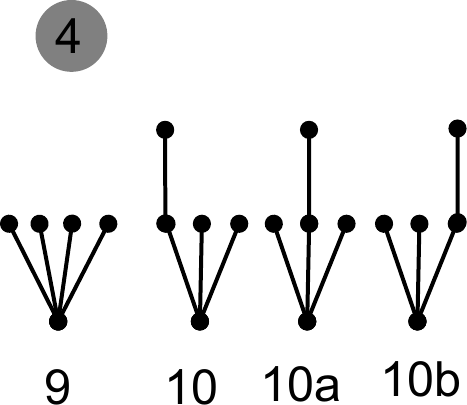} \ \
\includegraphics[width=0.3\linewidth]{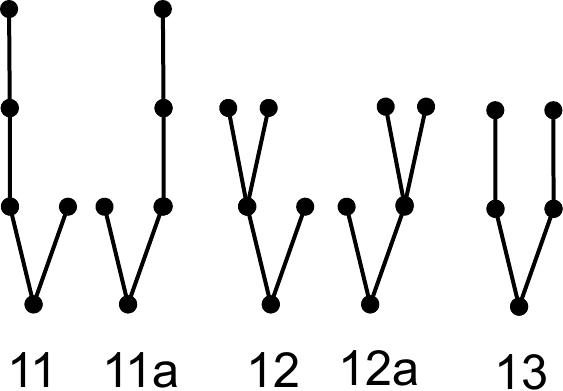} \ \
\includegraphics[width=0.27\linewidth]{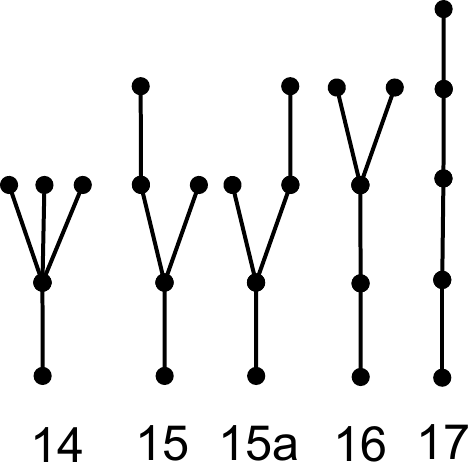}
}
\caption{rooted trees with 4 edges 
}
\label{4a}
\end{figure}

Root trees with 4 edges are shown in Fig. \ref{4a}.

Trees 4, 4a, 4b are isomorphic planar graphs, but not isomorphic rooted trees.

The same applies to pairs of graphs 11, 11a; 12, 12a; 15, 15a.

So, there are a total of 14 non-isomorphic rooted trees with 4 edges.

\subsection{Rooted trees with 5 edges}
\begin{figure}[ht!]
\center{\includegraphics[width=0.3\linewidth]{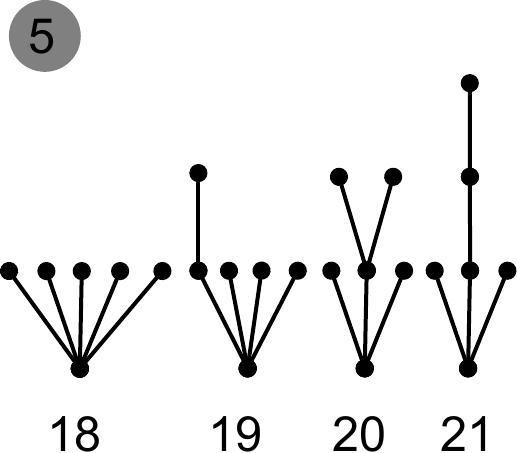} \ \
\includegraphics[width=0.3\linewidth]{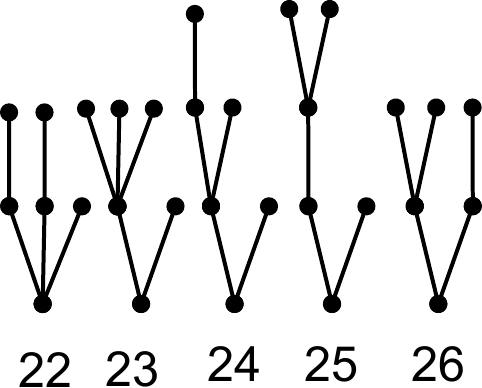}
}
\caption{rooted trees with 5 edges (part 1)
}
\label{5a}
\end{figure}

\begin{figure}[ht!]
\center{\includegraphics[width=0.3\linewidth]{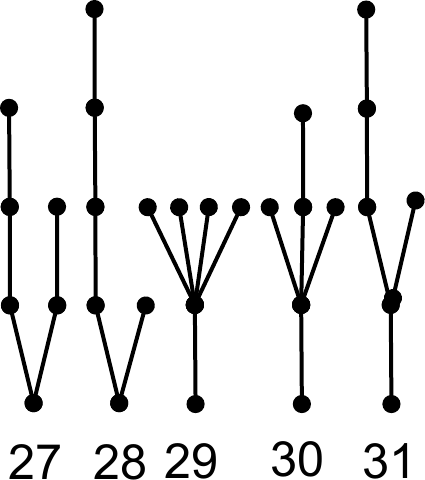} \ \
\includegraphics[width=0.3\linewidth]{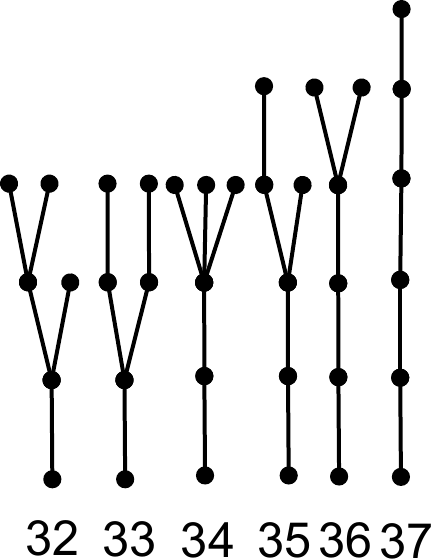}
}
\caption{rooted trees with 5 edges (part 2)
}
\label{5c}
\end{figure}

Root trees with 5 edges are shown in Fig. \ref{5a},\ref{5c}. Only one root tree in each plane tree class is given here.

For tree 19 there are 4 non-isomorphic rooted trees, for 20 there are three rooted trees, and for 21 there are also 3 rooted trees.

There are 3 non-isomorphic root trees for tree 22, 2 root trees for 23, 4 for 24, 2 for 25, and 2 root trees for 26.

For tree 27 there are 2 non-isomorphic root trees, for 28 - 2 root trees, and for 29 - one, for 30 - three, and for 31 - 2 root trees.

For tree 32, there are 2 non-isomorphic root trees, for 33 - one root tree, and for 34 - one, for 35 - two, and for 36 and 37 one root tree each.

So, there are a total of 51 non-isomorphic rooted trees with 5 edges.

\section{Planar trees}

Since rooted trees are derived from planar trees, if you select a vertex on them, planar trees with a small number of vertices can be specified by root tree numbers.

\subsection{Planar trees with no more than 6 vertices}

We will get the following list:
\begin{enumerate}

\item
A planar tree with one vertex -- 1.

\item
A planar tree with two vertices -- 2.

\item
A planar tree with three vertices -- 3=4.

\item
Two planar trees with four vertices -- 5=7, 6=6a=8.

\item
Three planar trees with five vertices -- 9=14, 10= 10a= 10b= 12= 12a= 15= 15a= 16, 11= 11a= 13= 17.

\item
Six planar trees with six vertices -- 18=29, 19=23=30=34, 20=32, 21 =25=26=31=36, 22=24=33=35, 27=28=37.
\end{enumerate}

\subsection{Planar trees with 7 vertices}
To obtain a list of planar trees, we will glue the root trees so that the roots are glued to a common vertex. At the same time, it should be taken into account that the same tree can be obtained by gluing different root trees. To get rid of such ambiguity, we will assume that the resulting vertex is central.
By central we will understand the vertex that is at the greatest distance from the vertices of valence 1 (there can be several central vertices)

\begin{figure}[ht!]
\center{\includegraphics[width=0.4\linewidth]{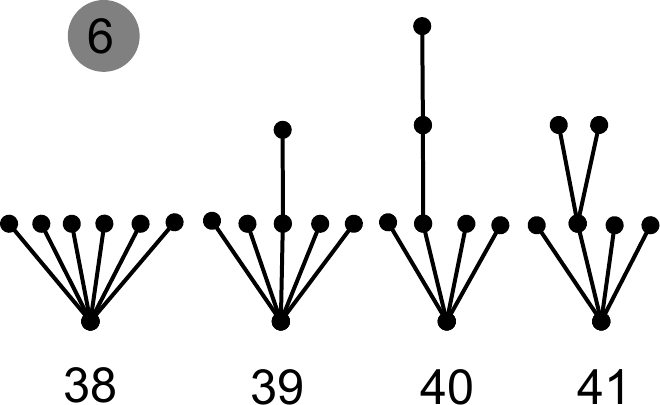} \ \
\includegraphics[width=0.4\linewidth]{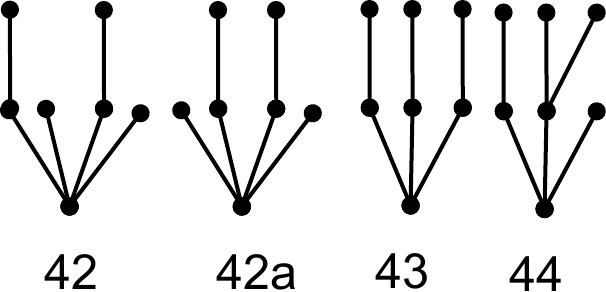}
}
\caption{planar trees with 7 vertices (part 1)
}
\label{t6b}
\end{figure}

\begin{figure}[ht!]
\center{\includegraphics[width=0.55\linewidth]{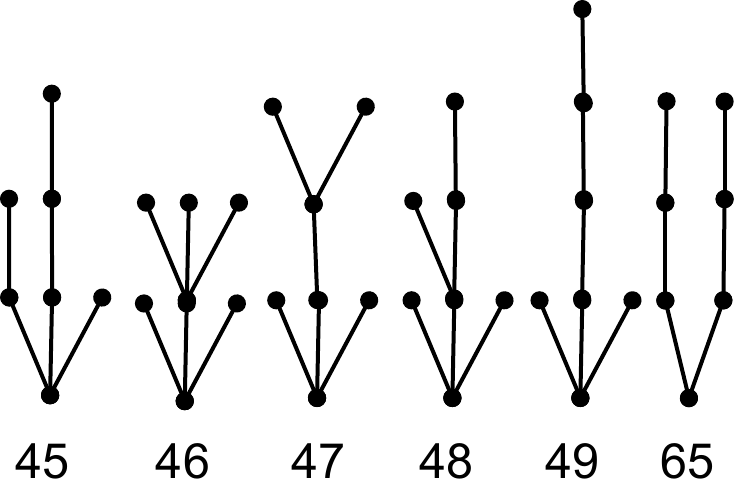}
}
\caption{planar trees with 7 vertices (part 2)
}
\label{t6c}
\end{figure}

All possible non-isomorphic planar trees with 7 vertices are shown in fig. \ref{t6b},\ref{t6c}.

So there are 14 planar trees with 7 vertices.

\subsection{Planar trees with 8 vertices}

We will match and list the trees in the same way as in the previous section.

\begin{figure}[ht!]
\center{\includegraphics[width=0.50\linewidth]{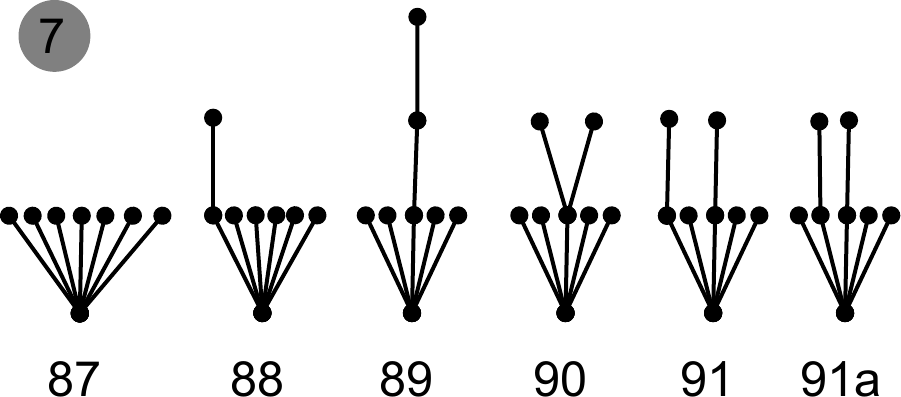} \ \
\includegraphics[width=0.35\linewidth]{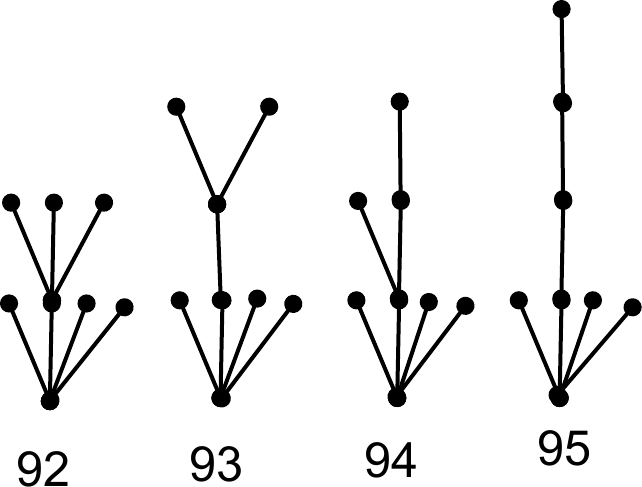}
}
\caption{planar trees with 6 vertices (part 1)
}
\label{t7b}
\end{figure}

\begin{figure}[ht!]
\center{\includegraphics[width=0.45\linewidth]{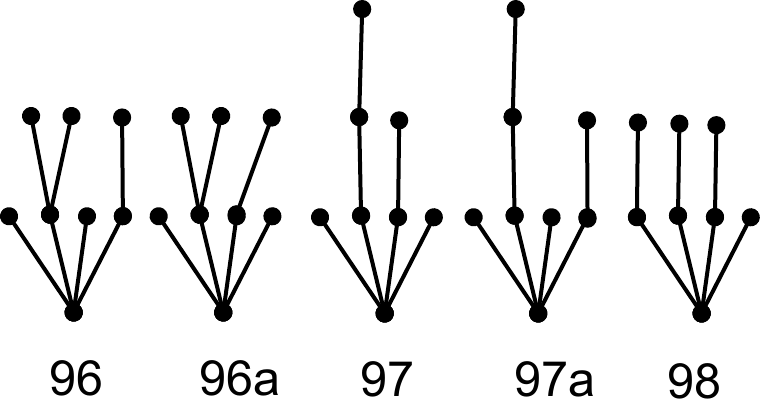}
\includegraphics[width=0.45\linewidth]{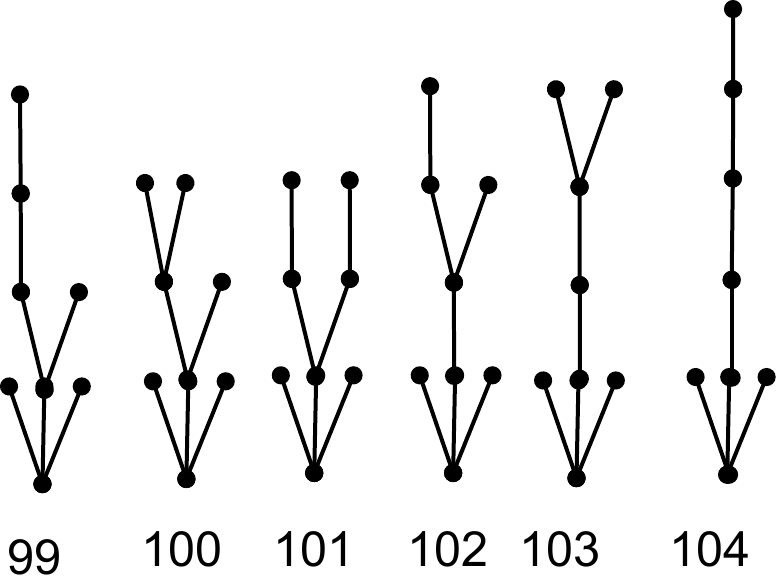}
}
\caption{planar trees with 6 vertices (part 2)
}
\label{t7d}
\end{figure}

\begin{figure}[ht!]
\center{\includegraphics[width=0.45\linewidth]{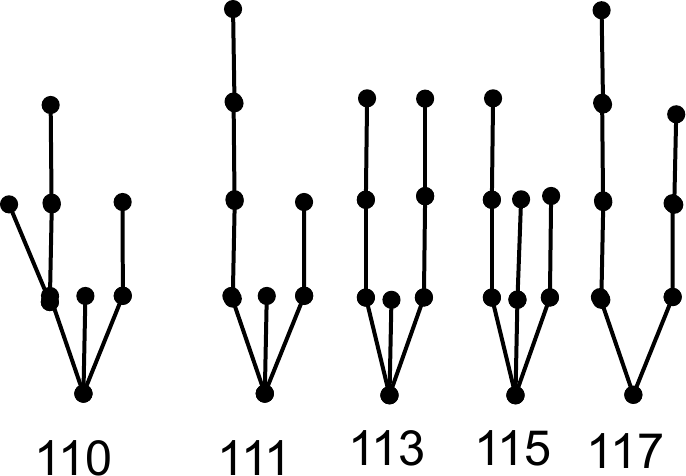}
}
\caption{planar trees with 6 vertices (part 3)
}
\label{t7e}
\end{figure}

All possible non-isomorphic planar trees with 8 vertices are shown in fig. 
\ref{t7b},
\ref{t7d},\ref{t7e}.

So there are 26 planar trees with 8 vertices.

\section{Morse flows on the 2-sphere}

Consider a Morse flow with one sink on the sphere. Such flows are gradient flows for the Morse function with one local maximum. Let us construct a graph $G$, the vertices of which are the sources and the edges are the stable manifolds of saddle points.

\begin{theorem} The graph $G$ constructed in this way is a complete topological invariant of the Morse flow. Moreover, the graph is a plane tree.
\end{theorem}

\textbf{\textbf{Proof.}} The fact that a graph is a complete topological invariant follows from \cite{Kybalko2018, Peixoto1973, prishlyak1997graphs, prishlyak2017morse}. Assume by contrast that the graph has a cycle. Then this cycle splits the sphere into two regions. Each of these areas has at least one drain. Then the total number of sinks is at least two, which contradicts the conditions of the theorem. The resulting contradiction completes the proof.

Note that since the Euler characteristic of the sphere is equal to 2, then, according to the Poincaré-Hopf theorem, the number of Morse flow sources with one sink is one more than the number of saddles.

\textbf{Consequence.} On the 2-, there are such a number of topologically non-equivalent Morse flows with one sink:

1 flow with one saddle,

1 flow with two saddles,

2 flows with three saddles,

3 flows with four saddles,

6 flows with five saddles,

14 flows with six saddles,

26 flows with seven saddles.

\section*{Conclusion}
\addcontentsline{toc}{chapter}{Conclusion}

In this paper, we have constructed a complete list of planar graphs with no more than 8 vertices. When constructing, we used rooted trees and checked them for isomorphism as planar graphs.


\end{document}